\documentclass[letter,11pt]{article}
\usepackage[utf8]{inputenc}
\usepackage[english]{babel}
\usepackage[all]{xy}
\usepackage{tikz}

\usepackage{amsmath}
\usepackage{amssymb}
\usepackage{amsthm}
\usepackage{mathrsfs}

\usepackage{mathtools}

\usepackage{amsfonts}

\usepackage[a4paper,top=2.5cm,bottom=2.5cm,left=2.5cm,right=2.5cm]{geometry}

\usepackage{color}

\renewcommand{\ker}{\operatorname{Ker}}

\newcommand{\del}{\partial}
\newcommand{\dbar}{\overline{\partial}}
\newcommand{\dbarb}{\overline{\partial}_b}

\newcommand{\su}{\mathfrak{su}}

\newcommand{\dint}{\ \!\mathrm{d}}
\newcommand{\dext}{\mathrm{d}}

\newcommand{\Z}{\mathbb{Z}}

\newcommand{\C}{\mathbb{C}}
\newcommand{\R}{\mathbb{R}}
\newcommand{\Q}{\mathbb{Q}}

\newtheorem{dfn}{Definition}
\newtheorem{exa}{Example}
\newtheorem{rmk}{Remark}
\newtheorem{lem}{Lemma}
\newtheorem{prp}{Proposition}
\newtheorem{thm}{Theorem}

\usetikzlibrary{matrix}

\begin{document}

\title{Levi-flat CR structures on compact Lie groups}

\author{Howard Jacobowitz \\ jacobowi@camden.rutgers.edu \and Max Reinhold Jahnke\footnote{During the development of this work, the second author received funding from FAPESP.} \\ jahnke@dm.ufscar.br}


\maketitle

\begin{abstract}
    Pittie \cite{pittie1988dolbeault} proved that the Dolbeault cohomology of all left-invariant complex structures on compact Lie groups can be computed by looking at the Dolbeault cohomology induced on a conveniently chosen maximal torus. We use the algebraic classification of left-invariant CR structures of maximal rank on compact Lie groups \cite{charbonnel2004classification} to generalize Pittie's result to left-invariant Levi-flat CR structures of maximal rank on compact Lie groups.
\end{abstract}

\section{Introduction}
\label{sec:introduction}
    
    In this work, we prove a version of the Leray-Hirsch theorem for CR principal bundles, which, combined with a result by Charbonnel-Kalgui \cite{charbonnel2004classification}, allows us to study the $\dbarb$ cohomology of left-invariant Levi-flat CR structures of maximal rank on compact Lie groups.
    
    Let $G$ be a connected and compact Lie group with Lie algebra $\mathfrak g$ and let $\C \mathfrak g$ be the complexification of $\mathfrak g$. We assume that $G$ is odd-dimensional and is endowed with a left-invariant Levi-flat CR structure $\mathcal V$ of CR codimension 1. Since $\mathcal V$ is left-invariant, there is a corresponding Lie algebra $\mathfrak h \subset \C \mathfrak g$ defined by the restriction of $\mathcal V$ to the identity of $G$. Notice that this is a one-to-one correspondence, that is, given a Lie subalgebra $\mathfrak h \subset \C \mathfrak g$, we define by left-translation an involutive vector bundle $\mathcal V \subset \C T G$ and if we assume that $\mathfrak h \cap \overline{\mathfrak h} = \{0\}$ the vector bundle $\mathcal V$ is an abstract CR structure. Notice that the rank of the vector bundle $\mathcal V$ is just the dimension of the Lie algebra $\mathfrak h$ and that the commutator bracket and the Lie algebra bracket agree. Further, because we use induced CR structure on a maximal torus, we restrict the rank of $G$ to be bigger than or equal to 3.
    
    We consider the $\dbarb$ complex on $G$
    $$ \dbarb : C^\infty(G; \Lambda^{0,q}) \to C^\infty(G; \Lambda^{0,q+1})$$ with $0 \leq q \leq \operatorname{rank} \mathcal V$ with cohomology spaces denoted by $H^{0,q}_{}(G; \mathcal V)$.

    Given $\xi \in \mathfrak h^\perp = \{u \in \C \mathfrak g^*: u(L) = 0~\forall L \in \mathfrak h\}$, $\xi \neq 0$, the \emph{Levi form} at  $\xi \in \mathfrak h^\perp$ is the Hermitian form on $\mathfrak h$ defined by $\mathcal L_{\xi}(L,M) = (1/2i)\xi([L,\overline{M}])$
    in which $L$ and $M$ are any smooth local sections of $\mathfrak h$. We shall say that $\mathcal L$ is \emph{Levi-flat} if given any point $\xi \in \mathfrak h^\perp$, with $ \xi \neq 0$, the Hermitian form $\mathcal L_{\xi}$ is zero. Notice that the Levi-form is left-invariant.
    
    By the classification theorem of Charbonnel-Kalgui \cite{charbonnel2004classification}, there is a maximal torus $T \subset G$ such that $\mathcal W_t \doteq \mathcal V_t \cap \C \operatorname{T}_t T$ defines a bi-invariant CR structure $\mathcal W = \bigcup_{t \in T} \mathcal W_t $ on $T$. We call $\mathcal W$ the toric component of $\mathcal V$ (relative to the maximal torus $T$). For each $g \in G$, we denote by $gT$ the set $\{gt: t \in T\}$ and by $g\mathcal W$ the pushforward by left-translation by $g$, that is, $g \mathcal W \doteq (L_g)_*(\mathcal W)$ with $L_g(x) = gx$. We usually denote the toric component $\mathcal W$ by its corresponding Lie algebra $\mathfrak m \subset \C \mathfrak t$ with $\mathfrak t$ the Lie algebra of $T$.
    
    We say that $\mathcal W$ (or $\mathfrak m$) satisfies the divisor condition {\bf(DC)} if there exist a basis $\{L_1, \ldots, L_n\}$ for $\mathcal W$ (or $\mathfrak{m}$) and constants $C, \rho > 0$ such that
    $$ \max_j |\widehat{L_j}(\xi)| \geq C(1 + |\xi|)^{-\rho}, \quad \forall \xi \in \Z^N,$$
    with $\hat L_j$ being the symbol of the vector field $L_j$ and $N$ the dimension of $G$.
    
    The following theorem is the main result of this work.     
    \begin{thm}
    \label{thm:main}
        Let $G$ be a connected, odd-dimensional, and compact Lie group endowed with a left-invariant Levi-flat CR structure $\mathcal V$ of maximal rank. Suppose that $\mathcal W$, the toric part of $\mathcal V$, satisfies the {\bf(DC)} condition, then there exist an isomorphism $H^{0,q}_{}(G; \mathcal V) \cong H^{0,q}_{}(T; \mathcal W)$.
    \end{thm}
    
    Notice that the toric part of left-invariant complex structures always satisfy the {\bf{(DC)}} condition and so our theorem is a generalization of a result by Pittie \cite{pittie1988dolbeault}.

    The proof of Theorem \ref{thm:main} is an adaptation to the CR case of a proof of the complex case we found in \cite{alexandrov2001vanishing}, and follows from the next two theorems: Theorem \ref{thm:main_leray_hirsch} is a version of Leray-Hirsch theorem for CR bundles assuming a cohomological condition and that the structure is Levi-flat; and Theorem \ref{thm:aux} shows that {\bf (DC)} condition  implies such a cohomological condition. 
    
    
    \begin{thm}
    \label{thm:main_leray_hirsch}
        Let $G$ be a connected and compact Lie group endowed with a left-invariant Levi-flat CR structure $\mathcal V$. A sufficient condition for the existence of an isomorphism for all $q$ $$H^{0,q}_{}(G; \mathcal V) \cong H^{0,q}_{}(T; \mathcal W)$$ is that for each $q$ and for every cohomology class $u$ in $H^{0,q}_{}(T; \mathcal W)$ there exist a cohomology extension $u'$ in $H^{0,q}_{}(G; \mathcal V)$ such that for all $g \in G$,  $u'$, when restricted to $gT$, defines a cohomology class $H^{0,q}_{}(gT; g\mathcal W)$.
    \end{thm}
    
    The proof of this theorem can be obtained using the general Leray-Hirsch theorem \cite[Theorem 9 of Section 7, Chapter 5]{spanier1994algebraic}. The hypothesis that $\mathfrak h$ is Levi-flat is used to show that the principal bundle $G$ with structure group $T$ and base space $G / T$ admits a CR trivialization. Instead of just citing that reference, we prove a version of Leray-Hirsch theorem for CR bundles using only the theory of partial differential equations and complex analysis.
    
    The next theorem shows that the {\bf(DC)} condition on the toric component of the CR structure guarantees that a cohomology extension exists.
    
    \begin{thm}
    \label{thm:aux}
        Let $G$ be a connected and compact Lie group endowed with left-invariant CR structure $\mathcal V$ of maximal rank and let $\mathcal W$ be its toric component. If $\mathcal W$ satisfies the {\bf(DC)} condition, then there exist a cohomology extension $H^{0,q}(T; \mathcal W) \to H^{0,q}(G; \mathcal  V)$ satisfying the conditions from Theorem \ref{thm:main_leray_hirsch}.
    \end{thm}
    
    The paper is organized in the following way. In Section \ref{sec:involutive_structures}, we briefly introduce the notation regarding locally integrable structures, the associated differential complexes, and their cohomology spaces.
    
    In Section \ref{sec:left_inv_cr}, we define the main object of study of this work, namely, the left-invariant CR structures of maximal rank. We also state the theorem of Charbonell-Kalgui \cite{charbonnel2004classification} that gives a classification of such CR structures in the language of Lie algebras. We finish this section with examples of such structures.
    
    In Section \ref{sec:lht-cr}, we briefly recall some concepts related to principal and CR bundles and prove a version of the Leray-Hirsch theorem for CR bundles.
    
    In Section \ref{sec:examples}, we discuss some examples to emphasize why some hypotheses are necessary and we also show that there are no Levi-flat CR structures of maximal rank on semisimple Lie groups.
    
    In Section \ref{sec:cr_struct_torus}, we show that {\bf (DC)} suffices to guarantee that cohomology class extensions exist. This sufficient condition is a Diophantine condition on the symbol of vector fields defining the toric part.
    
    In Section \ref{sec:op}, we describe some open problems that would nicely complement our results.

\section{Involutive structures}
\label{sec:involutive_structures}

Let $\Omega$ be a smooth and orientable manifold of dimension $N$. An involutive structure on $\Omega$ is a smooth subbundle $\mathcal V$ of the complexified tangent bundle $\C T \Omega$ of $\Omega$ such that the Lie bracket of any two smooth local sections of $\mathcal V$ is again a smooth section of $\mathcal V$. We denote the rank of $\mathcal V$ by $n$, and we denote by $\Lambda^k$ the bundle $\Lambda^k \C T^* \Omega$. If $W$ is a smooth vector bundle, we denote by $\mathcal C^\infty(\Omega; W)$ the space of sections of $W$ with smooth coefficients. 

For $0 \leq p \leq m \doteq N - n$, $0\leq q \leq n$ and each open set $U \subset \Omega$ there is a natural differential operator $\dext'$ associated to $\mathcal{V}$ defining a complex
\begin{equation}
    \label{h_complex}
    \dext': \mathcal C^\infty(U; \Lambda^{0,q}) \to \mathcal C^\infty(U; \Lambda^{0,q+1})
\end{equation}
with cohomology spaces denoted by $H^{0,q}_{} (U; \mathcal V).$ We refer to \cite{berhanu2008introduction} and \cite{treves1992hypo} for a detailed construction of this complex.

When $\mathcal V + \overline{\mathcal V} = \C T \Omega$, the structure $\mathcal V$ is called an \emph{elliptic structure}, when $\mathcal V \oplus \overline{\mathcal V} = \C T \Omega$, the structure $\mathcal V$ is called a \emph{complex structure} and the differential operator $\dext'$ is the usual $\dbar$ operator, and when $\mathcal V \cap \overline{\mathcal V} = \{0\}$, the structure $\mathcal V$ is called a \emph{(abstract) CR structure} and the differential operator $\dext'$ denoted by $\dbarb$. 




\section{Left-invariant involutive structures on compact Lie groups}
\label{sec:left_inv_cr}
    
    Let $G$ be a compact Lie group of dimension $N$. We denote by $\mathfrak g$ the Lie algebra of $G$. If $\mathfrak h \subset \mathfrak g$ is any subalgebra, we denote by $\C \mathfrak h$ its complexification. Let $L_g : G \to G$ be the left-multiplication by $g \in G$, that is $L_g(x) = g \cdot x$. Notice that the push-forward $(L_g)_*$ induces a vector bundle isomorphism between $G \times \C \mathfrak g$ and $\C T G$. Therefore we can identify subbundles of $\C T G$ with subbundles of $G  \times \C \mathfrak g$.
    
    A subbundle $\mathcal V \subset \C T G$ is said to be left-invariant if $(L_g)_* X \in \mathcal V_{g \cdot x}$ for all $X \in \mathcal V_x.$ We identify left-invariant involutive subbundles of $\C T G$ with Lie subalgebras of $\C \mathfrak g$. For a subalgebra $\mathfrak h \subset \C \mathfrak g$, we denote by $\mathcal V_\mathfrak h$ the associated subbundle. Notice that the rank of the bundle $\mathcal V_\mathfrak h$ is the complex dimension of $\mathfrak h$.
    
    \subsection{Left-invariant CR structures of maximal rank on compact Lie groups}

        In this section, we make a brief exposition of some results from \cite{charbonnel2004classification}. 
        
        Let $\mathfrak t \subset \mathfrak g$ be a maximal abelian subalgebra. We denote by $\Delta$ the set of roots of $\C \mathfrak t$ in $\C \mathfrak g$ and by $\Delta_+$ a maximal subset of positive roots of $\Delta$. For $\alpha \in \Delta$ we denote by $\mathfrak g_\alpha$ the eigenspace associated to $\alpha$. Each $\mathfrak{g}_\alpha$ is one dimensional \cite[Theorem 7.23]{hall2015lie}. We can easily show that
        \begin{equation}
        \label{eq:positive_roots_algebra}
            \mathfrak b_{\mathfrak t} = \bigoplus_{\alpha \in \Delta_+} \mathfrak g_\alpha
        \end{equation} is a Lie subalgebra of $\C \mathfrak g$.
        
        Since $\mathfrak{t}$ is abelian, any choice of vector space $\mathfrak{m} \subset \C \mathfrak{t}$ is a Lie algebra and $\mathfrak{h} = \mathfrak m \oplus \mathfrak{b}_\mathfrak{t}$ is a Lie algebra. If $\mathfrak{m}$ is elliptic (resp. CR), then $\mathfrak{h}$ is elliptic (resp. CR). Also, notice that $\mathfrak{b}_\mathfrak{t}$ is an ideal of $\mathfrak h$.
        
        Let us focus on the CR case. We denote by $d$ the dimension of $\mathfrak t$ and, since the dimension of each $\mathfrak g_\alpha$ is one, we can easily see that $N = d + 2l$ with $l$ being the number of elements of $\Delta_+$. If $\mathfrak m$ is a Lie subalgebra of $\C \mathfrak t$ of dimension $[d/2]$\footnote{For $x \in \R$, we define $[x] = \max \{ n \in Z: n \leq x \}$.} and such that its intersection with $\mathfrak g$ is null, the sum $\Phi(\mathfrak m)$ of $\mathfrak m$ and $\mathfrak b_{\mathfrak t}$ is a subalgebra of dimension $[N/2]$ and null intersection with $\mathfrak g$.
        
        We recall the two subalgebras $\mathfrak h_1$ and $\mathfrak h_2$ of $\mathfrak g$ are called conjugate if there exist a $g \in G$ such that $Ad(g) \mathfrak h_1 = \mathfrak h_2$. Here $\operatorname{Ad}(g)$ is the differential at $e \in G$ of the map $C_g : G \to G$ given by $C_g(h) = ghg^{-1}$. 
    
        \begin{dfn}
            We shall say that a subalgebra $\mathfrak h \subset \C \mathfrak g$ is of type CR0 relative to $\mathfrak g$ if it is conjugate to a subalgebra of the form
            \begin{equation}
            \label{eq:CR0}
                \Phi(\mathfrak m) = \mathfrak m \oplus \mathfrak b_{\mathfrak t} 
            \end{equation}
            in which $\mathfrak m$ is a subspace of dimension $[d/2]$ of $\C \mathfrak t$ and null intersection with $\mathfrak g$.
        \end{dfn}
    
        When $N$ is odd, we introduce the set $\mathcal D(\Delta_+)$ of all elements of the form $(\alpha, \mathfrak m, x, t)$ with $\alpha$ a simple root in $\Delta_+$, $\mathfrak m$ a subspace of dimension $[l/2]$ of the kernel of $\alpha$ with null intersection with $\mathfrak g$ and such that $x$ is an element not null of $\mathfrak g_\alpha$ and $t \in \mathfrak t$. For every $(\alpha, \mathfrak m, x, t)$ in  $\mathcal D(\Delta_+)$ we define
        \begin{equation}
        \label{eq:CR1}
            \Theta(\alpha, \mathfrak m, x, t) = \mathfrak m \oplus \bigoplus_{\beta \in \Delta_+ \backslash \{\alpha\}} \mathfrak g_\beta \oplus \operatorname{span}_\C\{t+x\}.
        \end{equation}
        
        Therefore $\Theta(\alpha, \mathfrak m, x, t)$ is a subalgebra of $\C \mathfrak g$ of dimension $[N/2]$ and null intersection with $\mathfrak g$.
    
    \begin{dfn}
        We shall say that a subalgebra $\mathfrak h \subset \C \mathfrak g$ is of type CR1 relative to $\mathfrak g$ if it is conjugate to a subalgebra of the form $\Theta(\alpha, \mathfrak m, x, t)$ in which $(\alpha, \mathfrak m, x, t)$ is an element of $\mathcal D(\Delta_+)$.
    \end{dfn}
    
    Without loss of generality, we take $\mathfrak h$ to be of the form \eqref{eq:CR0} or \eqref{eq:CR1} instead of conjugate to this form.
    
    For an algebra $\mathfrak{h}$ of type CR0 or CR1, we call the subalgebra $\mathfrak{m}$ the toric part of $\mathfrak{h}$.
    
    In Section \ref{sec:examples}, simple examples of both types of CR structures are given for the group $SU(2)$.
    
    The main theorem of \cite{charbonnel2004classification} is the following:
    
    \begin{thm}
        Let $\mathcal V$ be a left-invariant CR structure of maximum rank over a Lie group $G$. Then the rank of $\mathcal V$ is $[N/2]$. If $N$ is even, then the corresponding subalgebra $\mathfrak h$ given by the set of all left-invariant vector fields of $\mathcal V$ is a complex structure and it is of type CR0. If $N$ is odd, the corresponding subalgebra $\mathfrak h$ is of type CR0 or CR1.
    \end{thm}

\begin{thm}
A semisimple Lie group $G$  does not admit a Levi-flat CR structure of type CR0.
\end{thm}
\begin{proof}
Recall that we limit ourselves to CR structures of maximal rank, i.e., the codimension of $\mathcal V \oplus \bar{\mathcal V}$ is one.  Thus we only consider odd dimensional Lie groups.

Let $\mathfrak g$ be the complexified Lie algebra of the group.  If $\mathfrak g$ is not semisimple, then $G$ is not semisimple.  And to show $\mathfrak g$ is not semisimple 
 it suffices to show that $[\mathfrak g, \mathfrak g]$ is a proper subset of $\mathfrak g$.
 A CR structure of type CR0 has the form
\begin{equation*}
\mathcal V = \mathfrak m \oplus \bigoplus_{\beta \in \Delta^+} \mathfrak g_\beta
\end{equation*}
where $\mathfrak m \oplus \bar{\mathfrak m}$ is of codimension one in some maximal toral subalgebra of $\mathfrak g$.  Call this toral subalgebra $\mathfrak t$.  Note that there is an  element $T$ in the real Lie algebra of $G$ with 
\begin{equation*}
\mathfrak t = \{ T \} \oplus  \mathfrak m \oplus \bar {\mathfrak m} 
\end{equation*}

We have the corresponding decompositions
\begin{equation*}
\mathfrak g = \mathfrak t \oplus \bigoplus_{\beta \in \Delta} \mathfrak g_\beta = \{T\} \oplus \mathcal{V}\oplus \mathcal {\bar V}.
\end{equation*}
From
\[
[T, {\mathfrak m} ]= 0
\]
and
\[
[T, {\mathfrak g}_\beta ] \subset {\mathfrak g}_\beta
\]
we see that
\[
[ T, \mathcal V ] \subset \mathcal V  \mbox{  and  }  [T, \bar { \mathcal V } ]\subset \bar { \mathcal V }.
\]
The involutivity  and Levi flatness conditions give us
\[
[ \mathcal V \oplus \bar{\mathcal V} , \mathcal V \oplus \bar{\mathcal V} ] \subset  \mathcal V \oplus \bar{\mathcal V}
\]
and so we now have
\[
[\mathfrak g ,\mathfrak g]\subset  \mathcal V \oplus \bar{\mathcal V}
\]
which is a proper subset of $\mathfrak g$.
\end{proof}
Note that this proof does not apply to structures of type CR1.  In place of $[T,{\mathcal V}]\subset \mathcal V $ we only obtain
\[
[ T,{\mathcal V} ]\subset {\mathcal V}\oplus {\mathfrak g} _\alpha
\]
(in the notation of Definition 2), from which we cannot conclude
\[
[ T, \mathcal V\oplus \bar{\mathcal V} ] \subset  \mathcal V \oplus \bar{\mathcal V}.
\]
It is fortunate that the proof breaks down for CR1:
\begin{exa}
Consider the usual basis of $\mathfrak{sl_2}$
\[
T=\left(
\begin{array}{cc}
1&0\\
0&-1
\end{array}
\right)
,\ \ \ X=\left(\begin{array}{cc}
0&1\\0&0
\end{array}
\right) ,\ \ \ Y=\left(\begin{array}{cc}
0&0\\
1&0\end{array}
\right) 
\]
with bracket products
\[
[T,X]=2X,\ \ \ [T,Y]=-2Y,\ \ \ [X,Y]=T.
\]
Thus
\[
\bf{C}\otimes \mathfrak{sl_2}=\mathfrak{M}\oplus \mathfrak{g_2} \oplus \mathfrak{g_{-2}}
\]
with
\[
\mathfrak {m} =\{ T \},\ \ \  \mathfrak{g_2}=\{ X \} ,\ \ \ \mathfrak{g_{-2}}=\{ Y \}.
\]
Take $L=X+iT$ and note that $L$ is of type CR1.  Also note that
\[
[L,\bar{L} ]=2i(L+\bar{L}).
\]
So the semi-simple group $SL_2$ admits a Levi flat CR structure (of type CR1).

\end{exa}
This example is taken from Section 4 of \cite{bor2019left-invariant}.  It follows from the methods of that section that up to CR equivalence this is the only  left-invariant Levi-flat CR structure on $SL_2$.

\section{Leray-Hirsch theorem for CR bundles}
\label{sec:lht-cr}

Let $G$ be an odd-dimensional compact Lie group and let $\mathfrak h$ be a Lie algebra in the form \eqref{eq:CR0} or \eqref{eq:CR1}. We define $T = \exp_G(\mathfrak t)$ and thus $T$ is maximal torus endowed with a maximal CR structure $\mathfrak m$. We have that the inclusion $i : T \subset G$ is a map compatible with the involutive structures on $T$ and on $G$, meaning that $i_*(\mathfrak{m}) \subset \mathfrak{h}$. We also have that the quotient map $\pi : G \to \Omega \doteq G / T$ induces a complex structure on $\Omega$, that is $\Omega$ is endowed with the involutive structure $\pi_*(\mathfrak{h})$ (see \cite[Problem 2.57]{gadea2013analysis} for details). 

In the following, we recall the definition of principal bundle, and we prove that $T \xrightarrow{i} G \xrightarrow{\pi} \Omega$ is a smooth principal bundle. We also prove that if $\mathfrak h$ is a Levi-flat CR subalgebra 
then the local smooth trivialization for the bundle $G$ is compatible with the involutive structure on the fiber and on the base space.

\begin{dfn}
    A \emph{principal fiber bundle} $P$ with structure group $H$ is a triple $(P,H,R)$ with $P$ being a manifold, $H$ a Lie group and $R : P \times H \to P$  a smooth right action of $H$ on $P$ satisfying:
    \begin{enumerate}
        \item $\Omega \doteq P / H$ has a manifold structure making the projection $\pi : P \to \Omega$ smooth;
        \item $P$ is locally trivial, i.e., there is an open cover $\{ U_j\}$ of $\Omega$ and diffeomorphisms $\Phi_j : \pi^{-1}(U_j) \to U_j \times H$ satisfying
        $$ \Phi_j^{-1} (x, hg) = \Phi_j^{-1} (x, h)g $$
        for all $x \in U_j$ and all $g,h \in H$.
    \end{enumerate}
\end{dfn}


\begin{exa}
    Let $G$ be a compact Lie group and let $K \subset G$ be a closed subgroup. Then $G$ can be regarded as a principal bundle with structure group $K$ endowed with a right action $ R : G \times K \to G$. The quotient by this action defines a manifold $\Omega = \{ gK : g \in G\}$ with the projection $\pi : G \to \Omega$ being smooth.
    
    By \cite[Corollary 10.1.11]{hilgert2011structure} it is possible to find a covering $\{U_j\}$ of $\Omega$ and smooth sections $\sigma_j$ of the quotient map $\pi : G \to \Omega$ such that
    \begin{equation}
    \label{eq:hilgert_diffeo}
        (u, t) \in U_j \times K \mapsto \Psi_j(u,t) \doteq \sigma_j(u)t \in \sigma_j(U_j)K
    \end{equation} is a diffeomorphism onto an open subset of $G$ which we denote by $V_j$.
    
    The local trivialization follows from the fact that $\Omega$ has an open cover $\{U_j\}$ and local sections $\sigma_j : U_j \to G$ for the projection $\pi$ such that $\Psi_j : (u,k) \in  U_j \times K \mapsto \sigma_j(u)k \in \pi^{-1}(U_j)$ is a diffeomorphism. We define $\Phi_j = \Psi_j^{-1}$ and notice that $$\Phi_j^{-1} (x, hg) = \Psi_j (x, hg) = \sigma_j(x) hg = \Psi_j (x, h)g = \Phi_j^{-1} (x, h)g.$$
\end{exa}

We define the analogous principal fiber bundle for CR structures.

\begin{dfn}
    A \emph{CR principal bundle} $P$ with structure group $H$ is a triple $(P,H,R)$ with $P$ being a CR manifold, $H$ a Lie group and $R : P \times H \to P$  a free CR right action of $H$ on $P$ satisfying:
    \begin{enumerate}
        \item $\Omega \doteq P / H$ has a CR structure making the projection $\pi : P \to \Omega$ a CR map;
        \item $P$ is locally CR trivial, i.e., there is an open cover $\{ U_j\}$ of $\Omega$ and CR diffeomorphisms $\Phi_j : \pi^{-1}(U_j) \to U_j \times H$ satisfying
        $$ \Phi_j^{-1} (x, hg) = \Phi_j^{-1} (x, h)g $$
        for all $x \in U$ and all $g,h \in H$.
    \end{enumerate}
\end{dfn}

\begin{prp}
\label{prp:holomorphic_principal_toric_bundle}
    Let $G$ be a compact Lie group endowed with a left-invariant Levi-flat CR structure $\mathfrak h$ of maximal rank and let $\mathfrak{m} \subset \mathfrak{h}$ be its toric part. Let $\mathfrak{t}$ be the Lie algebra of the maximal torus $T$ such that $\mathfrak{m} \subset \C \mathfrak{t}$. Then $G$ is CR principal bundle with structure group $T$ and base space $\Omega = G / T$.
\end{prp}

In order to prove Proposition \ref{prp:holomorphic_principal_toric_bundle}, we need: a few results:

\begin{lem}
    \label{lem:right_multiplication}
    Let $G$ be a compact Lie group endowed with a left-invariant CR structure $\mathfrak h$ of maximal rank and let $\mathfrak{m} \subset \mathfrak{h}$ be its toric part and let $T$ be the associated maximal torus. Let $G \times T$ be endowed with the CR structure given by $\mathfrak h \oplus \mathfrak m$. 
    Then $\mu : (g,t) \in G \times T \mapsto gt \in G$ is a CR map.
\end{lem}

\begin{proof}
    Let $X \oplus Y \in \mathfrak h \oplus \mathfrak m$ and notice that $$(\mu_*)_{(a,b)} (X \oplus Y) = (R_b)_*(X) + (L_a)_*(Y).$$
    Let $\mathfrak h_g$ denote $(L_g)_* \mathfrak h$, so $(L_a)_*(Y) \in \mathfrak h_{ab}$, hence we only need to verify that $(R_b)_*(X) \in \mathfrak h_{ab}$. 
    
    We write $X = X' + \sum_{\alpha \in \Delta_+} X_\alpha$ with $X' \in \mathfrak m$ and $X_\alpha \in \mathfrak g_\alpha$. 
    Clearly, we have that $(R_b)_*(X') \in \mathfrak h_{ab}$, and by linearity, we just need to know what happens with $(R_b)_*(X_\alpha)$. 
    Let $W \in \C T_b T$ and recall that $[W,X_\alpha] = \alpha(W)X_\alpha$. Also, notice that $$[W_{ab}, (R_b)_*(X_\alpha)] = [(R_b)_*(W_a), (R_b)_*(X_\alpha)] = (R_b)_*[W_a, X_\alpha] = (R_b)_*\alpha(W)X_\alpha.$$ Therefore, $(R_b)_*(X_\alpha) \in \ker \left( \operatorname{ad}_W - \alpha(W) \operatorname{id}\right) = \mathfrak{g_\alpha}$, with $\operatorname{ad}_W(X) = [W,X]$. Now, since $\dim_\C \mathfrak g_\alpha = 1$, we have that there exists a $\lambda_{b,\alpha} \in \C$ such that $$(R_b)_*(X_\alpha) = \lambda_{b,\alpha} X_\alpha.$$
\end{proof}

\begin{lem}
\label{lem:section_for_levi_flat}
Let $M$ be a smooth manifold endowed with a complex or a Levi-flat CR structure $\mathcal V$ and $N$ be a smooth manifold endowed with a complex structure $\mathcal S $. Let $f : M \to N$ be a map such that $f_*(\mathcal{V}_x) = \mathcal{S}_{f(x)}$ for all $x \in M$. Then for every $y \in N$ there exist an open neighborhood $U$ of $y$ and a map $\sigma : U \to M$ such that $f \circ \sigma(x) = x$ for all $x \in U$ and $\sigma_*(\mathcal S _x) \subset \mathcal V_\sigma(x)$ for all $x \in U$.
\end{lem}


\begin{proof}
Let $y \in N$, $x \in f^{-1}(\{y\})$ and $V \subset M$ be an open neighborhood of $x$ with coordinates $\phi : V \to V' \times V'' \subset \C^\nu \times \R^{n-\nu}$ with $V' \subset \C^\nu$ and $V'' \subset \R^{n - \nu}.$ When $\mathcal V$ is complex we take $\nu = n$. We choose the coordinates $$\phi = (z_1, \ldots, z_\nu, t_1, \ldots, t_{n - \nu})$$ such that over $V$ the involutive structure $\mathcal V$ is generated by the vector fields $$\frac{ \del }{ \del \overline{z}_1}, \ldots, \frac{ \del }{ \del \overline{z}_\nu}.$$

Let $U \subset N$ be a neighborhood of $y$ with coordinates $\psi : U \to U' \subset \C^l$ written as $$ \psi = (w_1, \ldots, w_l)$$ such that over $U$ the involutive structure $\mathcal S $ is generated by the vector fields $$\frac{ \del }{ \del \overline{w}_1}, \ldots, \frac{ \del }{\del \overline{w}_l} .$$
Now, by using these coordinates, we write a local representation of $f$, namely $\tilde f : V' \to U'$, as $\tilde f = \psi \circ f \circ \phi^{-1}$. Notice that since $f_*(\mathcal{V}_x) = \mathcal{S}_{f(x)}$ we have that $$f_*(\mathcal{V}_x + \overline{\mathcal V_x}) = \mathcal{S}_{f(x)} + \overline{\mathcal{S}_{f(x)}} = \C T_{f(x)}N.$$

Let $J \tilde f (z,t)$ denote the Jacobian matrix of $\tilde f$. Since $\tilde f_*$ is surjective, by shrinking the open sets $V, V', U$ and $U'$ if necessary, and by rearranging the indices, we can assume that first $l \times l$ block of the Jacobian matrix of $\tilde f$ is invertible. We write $z' = (z_1, \ldots, z_l)$, $(z'_0,z''_0,t_0) = \phi(x)$, $w_0 = \psi(y)$, and we define the function $F(z', z'',t) = (\tilde f(z',z'',t_0),z'',t)$. Therefore, from the Inverse Function Theorem there exist $F^{-1}$ and this function is holomorphic in the $z',z''$ and smooth in $t$. Let $F'(z,t) = (F(z),t)$, notice that this function is invertible, and define $\phi' = F' \circ \phi$. Notice that with $\phi'$ as new coordinate on $V$ we have that the local representation of $f$ is the projection on $z'$.

Therefore we can define a holomorphic section for $f$ by defining $\sigma' (z',z'',t) = (z', z''_0, t_0)$ and then defining $\sigma : U \to M$ as $\sigma = \phi^{-1} \circ \tilde \sigma \circ \phi$. Thus we have that $f \circ \sigma = \operatorname{Id}_U$ and by construction, since $\tilde \sigma$, $\tilde \phi$ and $\psi$ are compatible with the structures $\mathcal V$ and $\mathcal S $, we have that $\sigma_*(w) \in V_{\sigma(x)}$ for all $w \in W_x$.
\end{proof}

\begin{proof}[Proof of Proposition \ref{prp:holomorphic_principal_toric_bundle}]

    We apply Lemma 2 with $M = G$ and $N = \Omega (=G/T)$ thus have for every small open set $U$ and every $t \in T$, the CR diffeomorphism
     \begin{equation} \label{eq:local_cr_diffeomorphism}
        \Phi : (u,t) \in U \times T \mapsto \sigma(u)t \in \sigma(U)T = \pi^{-1}(U).
    \end{equation}

    
    Let $X \oplus Y \in \pi_*(\mathfrak{h}|_u \oplus \mathfrak m)$. We want to prove that $(\Phi_*)_{(u,t)}(X \oplus Y) \in \mathfrak{h}_{\sigma(u)t}$. First, we compute $(\Phi_*)_{(u,t)}(0 \oplus Y)$. Let $\alpha, \beta : (-\epsilon,\epsilon) \to G$ be two smooth curves satisfying $\alpha(0) = \beta(0) = t$ and $\alpha'(0) + i \beta'(0) = Y$. Then
    $$
    \begin{aligned}
        (\Phi_*)_{(u,t)}(0 \oplus Y) & = \left. \frac{\dext}{\dext s} \right|_{s=0}
        \sigma(u)\alpha(s) +i \left. \frac{\dext}{\dext s} \right|_{s=0} \sigma(u)\beta(s) \\
        & = (L_{\sigma(u)})_*(\alpha'(0)) + i (L_{\sigma(u)})_*(\beta'(0)) \\
        & = (L_{\sigma(u)})_*(Y) \in \mathfrak{h}_{\sigma(u)t}
    \end{aligned}
    $$
    
    Now, we compute $(\Phi_*)_{(u,t)}(X \oplus 0)$. Let $\alpha, \beta : (-\epsilon,\epsilon) \to G$ be two smooth curves satisfying $\alpha(0) = \beta(0) = u$ and $\alpha'(0) + i \beta'(0) = X$. Then
    $$
    \begin{aligned}
        (\Phi_*)_{(u,t)}(X \oplus 0) & = \left. \frac{\dext}{\dext s} \right|_{s=0}
        \sigma(\alpha(s))t + i \left. \frac{\dext}{\dext s} \right|_{s=0} \sigma(\beta(t))t \\
        & = (R_t)_* \circ \sigma_* \alpha'(0) + i (R_t)_* \circ \sigma_* \beta'(0) \\
        & = (R_t)_* \circ \sigma_*(X)
    \end{aligned}
    $$
    
    We use that $R_t$ and $\sigma$ are CR maps, and the proof is complete.
\end{proof}

\begin{rmk}
    Notice that the existence of local sections for $\pi$ implies that  \eqref{eq:local_cr_diffeomorphism} is a CR diffeomorphism and since the set $U \times T$ is endowed with a Levi-flat structure, we conclude that the existence of local sections implies that the original structure is Levi-flat. Therefore, the Levi-flat condition from Proposition \ref{prp:holomorphic_principal_toric_bundle} can not be dropped.
\end{rmk}

In the following, we show an example of a left-invariant CR structure of type CR0 such that the map $\mu$, as defined on Lemma \ref{lem:right_multiplication}, is a CR map, but, since such structure is not Levi-flat, it fails to be a CR principal bundle with structure group $T$.

\begin{exa}
Let  $G$ be the compact Lie group
$$ \operatorname{SU}(2) \doteq \left\{ \left(  \begin{matrix} z_1 & -\overline{z_2} \\ z_2 & \overline{z_1} \end{matrix} \right) : z_1, z_2 \in \C, ~|z_1|^2 + |z_2|^2 = 1 \right\}$$ and let $\mathbb{T} \subset G$ given by
$$
    \mathbb{T} = \left\{ \begin{pmatrix}
e^{i \theta} & 0 \\
0            & e^{-i \theta} \\
\end{pmatrix}: \theta \in \R \right\}.$$

It is easy to prove that $\mathbb{T}$ is a maximal torus. Recall that the Lie algebra of $\operatorname{SU}(2)$, which we denote by $\mathfrak{su}(2)$, is generated by
$$ X = \left(  \begin{matrix} 0 & i \\ i & 0 \end{matrix} \right), \quad Y = \left(  \begin{matrix} 0 & -1 \\ 1 & 0 \end{matrix} \right), \quad  T = \left(  \begin{matrix} i & 0 \\ 0 & -i \end{matrix} \right).$$

The vector field $L = X \otimes 1 - Y \otimes i$ defines a left-invariant CR structure $\mathfrak{h} = \operatorname{span}_\C \{L\} \subset \C \mathfrak{g}$. For $t \in \mathbb T$, we can easily compute $(R_t)_*(L) = (R_t)_*(X) \otimes 1 - (R_t)_*(Y) \otimes i$. In fact, we write $$t = \begin{pmatrix}
e^{i \theta} & 0 \\ 0  & e^{-i \theta} \\ \end{pmatrix}$$ for a $t \in \R$ and we get
$$
    (R_t)_*(X) = Xt = \left(  \begin{matrix} 0 & i \\ i & 0 \end{matrix} \right)  \begin{pmatrix} e^{i \theta} & 0 \\ 0  & e^{-i \theta} \\ \end{pmatrix} = \begin{pmatrix} 0 & ie^{-i\theta} \\ ie^{i\theta} & 0 \end{pmatrix},
$$ and 
$$
    (R_t)_*(Y) = Yt = \left(  \begin{matrix} 0 & -1 \\ 1 & 0 \end{matrix} \right)  \begin{pmatrix} e^{i \theta} & 0 \\ 0  & e^{-i \theta} \\ \end{pmatrix} = \begin{pmatrix} 0 & -e^{-i\theta} \\ e^{i\theta} & 0 \end{pmatrix}
$$ thus $$(R_t)_*(L) = \begin{pmatrix} 0 & ie^{-i\theta} \\ ie^{i\theta} & 0 \end{pmatrix} \otimes 1 - \begin{pmatrix} 0 & -e^{-i\theta} \\ e^{i\theta} & 0 \end{pmatrix} \otimes i$$
and we easily verify that $(R_t)_*(L) = e^{-i\theta}L$. Thus $\mu$ is a CR map.
\end{exa}

We will use the following special case of the Künneth formula.

\begin{prp}
\label{lem:skf}
    Let $D \subset \C^n$ be a polydisc with its natural complex structure $\mathcal V$ and let $U$ be a smooth manifold endowed with a CR structure $\mathcal W$. We assume that there are globally defined linearly independent $\dbarb$-closed forms $\zeta_1, \ldots, \zeta_m$ spanning $\mathcal W^\perp$. Then $$H^{0,q}(D \times U; \mathcal V \oplus \mathcal W) \cong \mathcal O(D) \otimes H^{0,q}(U; \mathcal V)$$ with $\mathcal O(D)$ the set of all holomorphic functions on $D$.
\end{prp}

\begin{proof}
    We denote by $\dbar$ the differential operator associated to $\mathcal V$, by $\dbarb$ the differential operator associated to $\mathcal W$, and by $\dext'$ the differential operator associated to $\mathcal V \oplus \mathcal W$. Notice that $\dext' = \dbar + \dbarb$ and since $\dext' \circ \dext' = 0$ we conclude that $\dbar \circ \dbarb + \dbarb \circ \dbar = 0$. We denote by $\pi_D$ the projection on the first variable and by $\pi_U$ the projection on the second variable.

    Notice that if $f \in \mathcal O(D)$ and $[u] \in H^{0,q}(U; \mathcal V)$, then $\pi_D^*(f) \pi_U^*(u)$ defines an element in $H^{0,q}(D \times U; \mathcal V \oplus \mathcal W)$ and since the elements of the form $\pi_D^*(f) [\pi_U^*(u)]$ form a basis for $\mathcal O(D) \otimes H^{0,q}(U; \mathcal V)$ we have a homomorphism $$\psi : \mathcal O(D) \otimes H^{0,q}(U; \mathcal V) \mapsto H^{0,q}(D \times U; \mathcal V \oplus \mathcal W).$$ We will prove that $\psi$ is an isomorphism.
    
    First we prove that $\psi$ is surjective. Let $[u] \in H^{0,q}(D \times U; \mathcal V \otimes \mathcal W)$. We write $u = \sum_{|J| + |K| = q} u_{JK} \dint \overline{z}_J \wedge \tau_K$. We can decompose $u$ as a finite sum $u = \sum_{j+k = q} u_{jk}$ with $u_{jk} = \sum_{|J| = j, |K| = k} u_{IJ} \dint \overline{z}_J \wedge \tau_K$. Let $u_{jk}$ be such that $j$ is as big as possible with $u_{jk} \neq 0$.
    
    From $\dext' u = 0$ and the maximality of $j$ we have that $\dbar u_{jk} = 0$. Since $D$ is a polydisc, we can find $v_{jk}$ in $D$ with parameters in $U$ such that $\dbar v_{jk} = u_{jk}$. Now we can take $u' = u - \dbar v_{jk}$. Notice that $u'$ and $u$ are in the same cohomology class. Now, let $u'_{jk}$ be a nonzero form such that $j$ is as big as possible. This $j$ is less than the one obtained with $u_{jk}$. By repeating this process, after a finite number of steps we find a $\tilde u$ that can be represented as $\tilde u = \sum_{|K| = q} \tilde u_{K} \tau_K$ with $\tilde u_{K} \in C^\infty(D \times U)$ holomorphic in $D$.
    
    Now we write the Taylor series at $0 \in D$ for each $\tilde u_{K}$ and we obtain $$ \tilde u_{K} = \sum_{\alpha \in \Z_+^n} \frac{\del^\alpha u_{K}}{\del z^\alpha}(0,t) \frac{z^\alpha}{\alpha!}, \qquad z \in D, t \in U$$ with uniform convergence over compact sets. Notice that $$\tilde u = \sum_{\alpha \in \Z_+^n} \frac{z^\alpha}{\alpha!} \left(\sum_{|K|=q} \frac{\del^\alpha u_{K}}{\del z^\alpha}(0,t) \tau_K\right)$$ and each term $\sum_{|K|=q} \frac{\del^\alpha u_{K}}{\del z^\alpha}(0,t) \tau_K$ defines a cohomology class in $ H^{0,q}(U; \mathcal V)$. This proves surjectivity.
    
    Now assume that $[u] \in \mathcal O(D) \otimes H^{0,q}(U; \mathcal V)$ is such that $\psi([u]) = 0$. This means that there exist $v$ such that $\dext' v = u$. We write $v = \sum_{|J| + |K| = q - 1} v_{JK} \dint \overline z_I \wedge \tau_K$ and we see that $\dbar v_{JK} = 0$ if $|J| \geq 1$ and with a process similar to the one we did in the sujectivity case conclude that we have a solution $\tilde v$ with $\dbar \tilde v = u$ with $\tilde v \in \mathcal O(D) \otimes H^{0,q-1}(U; \mathcal V)$.
\end{proof}

Let $(P,G,R)$ be a principal CR bundle with base space $B$. Let $\mathcal U, \mathcal V$ and $\mathcal W$ be the CR structures of, respectively, $G,P$ and $B$ and denote by $r_x : H^{0,q}(P; \mathcal V) \to H^{0,q}(P_x; \mathcal V|_{P_x})$ the restriction to the fiber $P_x$. A homomorphism $e : H^{0,q}(G; \mathcal U) \to H^{0,q}(P; \mathcal V)$ is called a \emph{cohomology extension of the fiber} if for every $x \in P$ the composition
$H^{0,q}(G; \mathcal U) \xrightarrow{e} H^{0,q}(P; \mathcal V) \xrightarrow{r_x} H^{0,q}(P_x;\mathcal V|_{P_x})$ is an isomorphism. We also assume that if $V \subset P$ is a small neighborhood such that there exist a trivialization $\Psi : V \to U \times G$ with $U = \pi(V)$ and if $\operatorname{pr_G}$ is the projection $\operatorname{pr_G} : U \times G \to G$, then
\begin{equation}
\label{equation:should_be_a_lemma}
    \Psi^* \circ \operatorname{pr_G}^* (u) =  e(u)
\end{equation}
    for all $u \in H^{0,q}(G; \mathcal U)$ and any $q = 0, 1, 2, \ldots.$

\begin{thm}[Leray-Hirsch theorem for CR principal bundles]
\label{thm:lhcr}
    Let $(P,H,R)$ be a principal CR bundle. Let $\mathcal V$ be the CR structure of $P$ and assume that there are globally defined linearly independent $\dbarb$-closed forms $\zeta_1, \ldots, \zeta_m$ spanning $\mathcal V^\perp$. Suppose that each $H^{0,q}(H)$ is finite-dimensional and that there is a cohomology extension $e : H^{0,q}(H) \to H^{0,q}(P)$, then there is an isomorphism
    $$ H^{0,q}_{}(P; \mathcal V) \cong
        \sum_{r+s=q} H^{0,s}_{}(H; \mathcal U) \otimes H^{0,r}_{}(\Omega; \mathcal W),$$
    with $\mathcal U$ the CR structure of $H$ and $\mathcal S$ the CR structure of the base space $\Omega$.
\end{thm}

\begin{proof} Let $\{U_j\}$ be an open covering of $\Omega$ by sets such that each $U_j$ is biholomorphic to a polydisc and there exists a local CR trivialization $\Phi_j : V_j \to U_j \times \Omega$ of $\Omega$ with $V_j = \pi^{-1}(U_j)$. Notice that $\{V_j\}$ is an open covering for $P$. Now we have an isomorphism
$$H^{0,q}_{}(V_j; \mathcal V) = H^{0,q}_{}(U_j \times H; \mathcal S \oplus \mathcal U) $$
and we have a homomorphism
$$ \psi : \sum_{r + s = q} H^{0,r}_{}(U_j; \mathcal W) \otimes H^{0,s}_{}(H; \mathcal U) \to H^{0,q}_{}(U_j \times H; \mathcal S \oplus \mathcal U)$$ given by $[u] \otimes [v] \mapsto \pi_{U_j}^*(u) \wedge \pi_{H}^*(v)$. By Lemma \ref{lem:skf} this last homomorphism is an isomorphism. These two isomorphism, combined with the fact that $e$ is a cohomology extension, implies that we have for each $V_j$ an isomorphism $$\phi : \sum_{r + s = q} H^{0,r}_{}(U_j; \mathcal W) \otimes H^{0,s}_{}(H; \mathcal U) \to H^{0,q}_{}(U_j \times H; \mathcal S \oplus \mathcal U)$$ given by  $[u] \otimes [v] \mapsto \pi^*(u) \wedge e(v)$.

On the other hand, for any two sets $U, U' \in \{U_j\}$ we can use Mayer-Vietoris sequence to get an exact sequence
$$ \cdots \to H^{0,r}_{}(U \cup U'; \mathcal W) \to H^{0,r}_{}(U; \mathcal W) \oplus H^{0,r}_{}(U'; \mathcal W) \to H^{0,r}_{}(U \cap U'; \mathcal W) \to H^{0,r+1}_{}(U \cup U'; \mathcal W) \to \cdots.$$

Now, we tensor each term with the finite vector space $H^{0,s}_{}(H; \mathcal U)$ and we sum all terms satisfying $r + s = q$ and we obtain another exact sequence. We define $V_j = \pi^{-1}(U_j)$. We take $V = \pi^{-1}(U)$ and $V' = \pi^{-1}(U')$. Now we have two sequences, which we show side by side emphasizing the isomorphism between each element. Notice that the only missing isomorphism is the term in the middle.

\begin{tikzpicture}
  \matrix (m) [matrix of math nodes,row sep=2em,column sep=0.5em,minimum width=1em]
  {
    \vdots & \vdots \\
    H^{0,q-1}_{}(V; \mathcal V) \oplus H^{0,q-1}_{}(V'; \mathcal V) &
        \sum H^{0,s}_{}(H; \mathcal U) \otimes H^{0,r}_{}(U; \mathcal W) \oplus \sum H^{0,s}_{}(H; \mathcal U) \otimes H^{0,r}_{}(U'; \mathcal W)  \\
    H^{0,q-1}_{}(V \cap V'; \mathcal V) &
        \sum H^{0,s}_{}(H; \mathcal U) \otimes H^{0,r}_{}(U \cap U'; \mathcal W)  \\
     H^{0,q}_{}(V \cup V'; \mathcal V) &
        \sum H^{0,s}_{}(H; \mathcal U) \otimes H^{0,r}_{}(U \cup U'; \mathcal W)  \\
     H^{0,q}_{}(V; \mathcal V) \oplus H^{0,q}_{}(V'; \mathcal V) &
        \sum H^{0,s}_{}(H; \mathcal U) \otimes H^{0,r}_{}(U; \mathcal W) \oplus \sum H^{0,s}_{}(H; \mathcal U) \otimes H^{0,r}_{}(U'; \mathcal W) \\
    H^{0,q}_{}(V \cap V'; \mathcal V) &
        \sum H^{0,s}_{}(H; \mathcal U) \otimes H^{0,r}_{}(U \cap U'; \mathcal W)  \\
     \vdots & \vdots \\};
  \path[-stealth]
    (m-1-1) edge node {} (m-2-1)
    (m-2-1) edge node {} (m-3-1)
    (m-3-1) edge node {} (m-4-1)
    (m-4-1) edge node {} (m-5-1)
    (m-5-1) edge node {} (m-6-1)
    (m-6-1) edge node {} (m-7-1)
    (m-1-2) edge node {} (m-2-2)
    (m-2-2) edge node {} (m-3-2)
    (m-3-2) edge node {} (m-4-2)
    (m-4-2) edge node {} (m-5-2)
    (m-5-2) edge node {} (m-6-2)
    (m-6-2) edge node {} (m-7-2)
    (m-2-1) edge node {} (m-2-2)
    (m-3-1) edge node {} (m-3-2)
    (m-4-1) edge [dashed,-] (m-4-2)
    (m-5-1) edge node {} (m-5-2)
    (m-6-1) edge node {} (m-6-2);
\end{tikzpicture}

Now, since this diagram is commutative, we can apply the Five Lemma and construct the missing isomorphism. Since $\Omega$ is compact, it has a finite covering, thus we can use a induction on the covering and construct an isomorphism defined in the union of all $\{U_j\}$ and $\{V_j\}$. That is, we obtain
$$ H^{0,q}_{}(P; \mathcal V) =
        \sum_{r+s=1} H^{0,s}_{}(H; \mathcal U) \otimes H^{0,r}_{}(\Omega; \mathcal W).$$
\end{proof}

\begin{lem}
    \label{lem:cohomology_extension}
    Let $G$ be a compact Lie group endowed with a left-invariant involutive structure $\mathfrak h$ of the form $$\mathfrak h = \mathfrak{m} \oplus \mathfrak u $$ with $\mathfrak m \subset \C \mathfrak t$ and $\mathfrak{t}$ the Lie algebra of a maximal torus $T$ and $\mathfrak u$ an ideal of $\mathfrak h$. Let $u \in C^\infty(T; \Lambda^{0,q})$ be a left-invariant and $\dbarb'$-closed, with $\dbarb'$ being the differential operator associated to $\mathfrak m$ on $T$. Then $u$ can be extended to a $\dbarb$-closed form in $G$ which restricts to a left-invariant form on each leaf $gT$.
\end{lem}

\begin{proof}
 Since $u$ is left-invariant we can regard it as an element of the dual of $\mathfrak m$ and we can extend $u$ as an element of $\C \mathfrak g^*$ by defining it as zero if any of its arguments are in $\mathfrak u$ or $\overline{\mathfrak u}$. We denote $\mathfrak u + \overline{\mathfrak u}$ as $\mathfrak m^\perp$. Let $w_g = (L_{g^{-1}})^* u$ and notice that $$(L_g)^* : H^{0,q}_{} (gT; \mathfrak m) \to H^{0,q}_{}(T; \mathfrak m )$$ is an isomorphism and so $w$ is a smooth $(0,q)$-form in $G$ which restricts to a cohomology class in each leaf $gT$.
    
    We are going to see that $w$ also is $\dbarb'$-closed. We just need to show that $w$ satisfies $\dext u(X_1, \dots, X_{q+1}) = 0$ if all the arguments are in $\mathfrak h$. Let $X_1, \ldots, X_{q+1} \in \mathfrak h$. Since the exterior derivative commutes with the pullback we have
    $$ \dext w_g = \dext [(L_{g^{-1}})^* u] = (L_{g^{-1}})^* (\dext u)$$
    we just need to do the following computation
    $$
    \begin{aligned}
    \dext u  (X_1,  \ldots, & X_{q+1})  = \sum_{j=1}^{q+1} (-1)^{j+1} X_j u (X_1, \ldots, \hat {X_j}, \ldots, X_{q+1}) \\
        & + \sum_{j < k} (-1)^{j+k+1} u ([X_j, X_k], X_1, \ldots, \hat {X_j}, \ldots, \hat {X_k}, \ldots, X_{q+1}) \\
        & = \sum_{j < k} (-1)^{j+k+1} u([X_j, X_k], X_1, \ldots, \hat {X_j}, \ldots, \hat {X_k}, \ldots, X_{q+1}). \\
    \end{aligned}
    $$
    
    Notice that by left-invariance, we have that $X_j u (X_1, \ldots, \hat {X_j}, \ldots, X_{q+1})$ is zero. If we assume that $X_j \in \mathfrak m$ for all $j$ then $(\dext  w_g)  (X_1,  \ldots, X_{q+1}) = 0$ because $\dbarb' u = 0$ and if we assume that $X_1 \in \mathfrak u$ and that $X_j \in \mathfrak u$ for some $j > 1$. Then
    $$
    \begin{aligned} 
    \sum_{j < k} (-1)^{j+k+1} & u([X_j, X_k], X_1, \ldots, \hat {X_j}, \ldots, \hat {X_k}, \ldots, X_{q+1}) \\ &= \sum_{1 < k} (-1)^{1+k+1} u([X_1, X_k], X_2, \ldots, \hat {X_j}, \ldots, \hat {X_k}, \ldots, X_{q+1}) = 0
    \end{aligned}
    $$
    because $\mathfrak u$ is an ideal and so $[X_1, X_k] \in \mathfrak u$.
    
    If we assume that two or more elements are $X_1,X_2 \in \mathfrak m^\perp$ we use the fact that $\mathfrak u$ is an ideal and an argument similar to the one above proves that we have $$(\dext  w_g)  (X_1(g),  \ldots, X_{q+1}(g)) = 0.$$ 
\end{proof}

\begin{rmk}
    Notice that if $\mathfrak h$ is a CR algebra of type CR0 or CR1, then we can use Lemma \ref{lem:cohomology_extension} with $\mathfrak u = \bigoplus_{\alpha \in \Delta_+} \mathfrak g_\alpha$ or $\mathfrak u = \bigoplus_{\beta \in \Delta_+ \backslash \{\alpha\}} \mathfrak g_\beta \oplus \operatorname{span}_\C\{t+x\}$. We refer to \cite{charbonnel2004classification} for a proof that, in both cases, $\mathfrak u$ is an ideal of $\mathfrak h$.
\end{rmk}

\begin{lem}
    \label{lem:cohomology_extension_inclusion}
    Let $G$, $T$, $\mathfrak h$ and $\mathfrak m$ be as in Lemma \ref{lem:cohomology_extension} and assume that every cohomology class of $H^{0,q}_{}(T; \mathfrak m)$ admits a left-invariant representative. Then the cohomology extension $$e : H^{0,q}_{}(T; \mathfrak m) \xhookrightarrow{} H^{0,q}_{}(G; \mathfrak h)$$ satisfies
    $$ \operatorname{pr}_T^*(u) = \Phi^*_U(e (u))$$
    for all cohomology classes $[u] \in H^{0,q}_{}(T; \mathfrak m)$ with $\Phi_U : (u, t) \in U \times T \mapsto \sigma(u)t \in \sigma(U)T$ being the local trivialization map and $\operatorname{pr}_T$ being the projection of $U \times T$ onto $T$.
\end{lem}

\begin{proof}
    
    Let $[u]$ be a cohomology class in $H^{0,q}(T; \mathfrak m)$ with $u$ left-invariant. Let $X_1 \oplus Y_1, \ldots, X_q \oplus Y_q \in \pi_*(\mathfrak h) \oplus \mathfrak m$.

    On the left hand side, we have $$\operatorname{pr}_T^*(u)(X_1 \oplus Y_1, \ldots, X_q \oplus Y_q) = u((\operatorname{pr}_T)_*(X_1 \oplus Y_1), \ldots, (\operatorname{pr}_T)_*(X_q \oplus Y_q)) = u(Y_1, \ldots, Y_q).$$

    And, on the right hand side, we have,
    $$
    \begin{aligned}
        \Phi^*_U(e (u))(X_1 \oplus Y_1, \ldots, X_q \oplus Y_q) = e(u)((\Phi_U)_*(X_1 \oplus Y_1), \ldots, (\Phi_U)_*(X_q \oplus Y_q)).
    \end{aligned}
    $$

Notice that $(\Phi_U)_*(X_j \oplus Y_j) = (R_t)_* \circ \sigma_*(X_j) + (L_{\sigma(u)})_*(Y_j)$ (see proof of Lemma \ref{prp:holomorphic_principal_toric_bundle}) and if we prove that $(R_t)_* \circ \sigma_*(X_j) \in \mathfrak m^\perp$ then $e(u)((\Phi_U)_*(X_1 \oplus Y_1), \ldots, (\Phi_U)_*(X_q \oplus 0)(\Phi_U)_*(X_q \oplus Y_q)) = 0$ and by linearity $e(u)((\Phi_U)_*(X_1 \oplus Y_1), \ldots, (\Phi_U)_*(X_q \oplus Y_q)) = u(Y_1, \ldots, Y_q)$.

Let $L_j \in \bigoplus_{\alpha \in \Delta} \mathfrak g_\alpha$ be such that $\pi_*(L_j) = X_j$. We can write $L_j = \sum_{\alpha \in \Delta_+} L_{j,\alpha}$ with $L_{j,\alpha} \in \mathfrak g_\alpha$ and $X_j = \sum_{\alpha \in \Delta_+} X_{j,\alpha}$ with $X_{j,\alpha} = \pi_*(L_{j,\alpha})$.

Notice that $\pi_*((R_t)_* \circ \sigma_*(X_{j,\alpha}) - L_{j,\alpha}) = \pi_*(\sigma_*(X_{j,\alpha})) - \pi_*(L_{j,\alpha}) = 0$ and thus $$\pi_*((R_t)_* \circ \sigma_*(X_{j,\alpha}) - L_{j,\alpha}) \in \C T_{\sigma(u)t} (\sigma(u)T) = \C \mathfrak t_{\sigma(u)t}.$$

Now, for any $W \in \C \mathfrak t_{\sigma(u)t}$, we have $$0 = [W, \pi_*((R_t)_* \circ \sigma_*(X_{j,\alpha}) - L_{j,\alpha}] = [W, \pi_*((R_t)_* \circ \sigma_*(X_{j,\alpha})] - [W,L_{j,\alpha}]$$ and so $[W, \pi_*((R_t)_* \circ \sigma_*(X_{j,\alpha})] = \alpha(W) L_{j,\alpha}$ for all $W$. This means that $\pi_*((R_t)_* \circ \sigma_*(X_{j,\alpha}) \in \mathfrak g_\alpha$ and thus $(R_t)_* \circ \sigma_*(X_j) \in \bigoplus_{\Delta_+} \mathfrak g_\alpha$. \end{proof}

Now we have all we need to prove Theorem \ref{thm:main_leray_hirsch}. Let $G$ be a connected and compact Lie group endowed with a left-invariant Levi-flat CR structure $\mathcal V$ which we represent by the Lie algebra $\mathfrak h$. We are assuming that that for every degree $q$ and every cohomology class $u$ in $H^{0,q}_{}(T; \mathfrak m)$ there exist a cohomology extension $u'$ in $H^{0,q}_{}(G; \mathcal V)$ such that, when restricted to $gT$, it defines a cohomology class $H^{0,q}_{}(gT; g\mathcal W)$. 

By Theorem \ref{thm:lhcr}, we have
$$ H^{0,q}_{}(G; \mathcal V) =
        \sum_{r+s=q} H^{0,s}_{}(T; \mathcal U) \otimes H^{0,r}_{}(\Omega; \mathcal W).$$
Notice that, from Chapter 8 of \cite{besse2008einstein}, $\Omega = G/T$ has positive first Chern class, so we can apply Corollary 11.25 of \cite{besse2008einstein} to obtain that $ H^{0,r}_{}(\Omega; \mathcal W) = 0$ for all $r > 0$. We conclude that 
$$ H^{0,q}_{}(G; \mathcal V) =
        H^{0,q}_{}(T; \mathcal U).$$



Let $G$ be a compact Lie group with Lie algebra $\mathfrak g$ and let $\mathfrak h \subset \C \mathfrak g$ be a Levi-flat CR algebra of maximal rank. That is, we are assuming that $2\dim_\C \mathfrak h + 1 = \dim_\C \mathfrak g$ and since $\mathfrak h$ is Levi-flat, we have that $\mathfrak h + \overline{\mathfrak h}$ is a Lie algebra so $\mathfrak k = (\mathfrak h + \overline{\mathfrak h}) \cap \mathfrak g$ is a real Lie algebra.
We can easily verify that $\mathfrak h + \overline{\mathfrak h} = \C \mathfrak k$ and the group $K = \exp_G(\mathfrak k)$ has a complex structure induced by $\mathfrak h$.

\begin{prp}
    Let $G$ be a compact Lie group with Lie algebra $\mathfrak g$. Let $\mathfrak h \subset \C \mathfrak g$ be a Levi-flat CR algebra of maximal rank and let $\mathfrak k = (\mathfrak h + \overline{\mathfrak h}) \cap \mathfrak g$ and suppose that $K = \exp_G(\mathfrak k)$ is closed. Then $$ H^{0,q}_{}(G; \mathfrak{h}) = H^{0,q}_{}(K; \mathfrak{h}) \otimes C^\infty(G/K).$$
\end{prp}

\begin{proof} The proof can be obtained by adapting the ideas from Theorem \ref{thm:main}.
\end{proof}

\begin{rmk}
    Let $G$, $K$ and $\mathfrak h$ be as in the proposition above. Notice that $\mathfrak h$ is a left-invariant complex structure and so it is of type CR0, therefore we can decompose $\mathfrak h$ as $\mathfrak h = \mathfrak m \oplus b_{\mathfrak t}$ with $\mathfrak t \subset \mathfrak k$ a maximal abelian subalgebra defining a maximal torus $T$. We have that $$H^{0,q}_{}(K; \mathfrak{h}) = H^{0,q}_{}(T; \mathfrak{m})$$ and so
    $$ H^{0,q}_{}(G; \mathfrak{h}) = H^{0,q}_{}(K; \mathfrak{h}) \otimes C^\infty(G/K) = H^{0,q}_{}(T; \mathfrak{m}) \otimes C^\infty(G/K).$$
\end{rmk}

\section{Examples}
\label{sec:examples}

        Recall that $\operatorname{SU}(2)$ is the group of all $2 \times 2$ unitary matrices having determinant 1. The Lie algebra of $\operatorname{SU}(2)$ is denoted by $\mathfrak{su}(2)$ and is the set of all skew-Hermitian and traceless $2 \times 2$ matrices. Let $\mathfrak k \subset \mathfrak{su}(2)$ be any one-dimensional subalgebra. Since every connected subgroup of $\operatorname{SU}(2)$ is closed and the rank of $\operatorname{SU}(2)$ is 1, we have that $K = \exp_{\operatorname{SU}(2)}(\mathfrak{k})$ is a maximal torus. Therefore, we can decompose $\C \mathfrak{su}(2)$ by using the roots associated to $\C \mathfrak{k}$. That is, we have
        \begin{equation*}
            \label{eq:su2_decomposition}
            \C \mathfrak{su}(2) = \C \mathfrak{k} \oplus \mathfrak{su}(2)_\alpha \oplus \mathfrak{su}(2)_{-\alpha}.
        \end{equation*}
    
        \begin{exa}
        \label{exa:cr1_su(2)}
            Let $X \in \mathfrak{k}$ be any element and let $Z \in \mathfrak{su}(2)_\beta$ be non-zero. We define $\mathfrak{h} = \operatorname{span}_\C \{X + Z\}$. If $X = 0$, the subalgebra $\mathfrak{h} = \mathfrak{su}(2)_\beta$ is of type CR0 and if $X \neq 0$, $\mathfrak{h}$ is of type CR1. These examples are not Levi-flat.
        \end{exa}

    The following example is actually a family of Levi-flat examples that will be useful.

    \begin{exa}
        \label{exa:maximum_rank_CR_structure_on_the_3-torus}
        Let $\mathbb{T} = \mathbb{T}^3$ be the 3-torus with Lie algebra $\mathfrak{t}$, with coordinates $(x,y,t)$, and endowed with the CR structure $\mathfrak{v} = \operatorname{span}_\C \{ L \}$ with $L = \del_x + \lambda \del_y + i \del_t$ and $\lambda \in \R$. We have that $\{L,\overline{L},\del_y\}$ is a basis for $\C \mathfrak{t}$. By taking $v = (1/2) (\dext x - i \dext t)$ and $w = -\lambda \dext x + \dext y$, we have a dual basis $\{v, \overline{v}, w\}$ for $\C \mathfrak{t}^*$. In this case, we have that $\Lambda^{0,1} = \operatorname{span}_\C \{ \dext v \}$, $C^\infty(\mathbb{T}; \Lambda^{0,1}) \cong \{f v\}$ and $\dbarb f = (Lf) u$ for all $f \in C^\infty(\mathbb{T})$. Notice that if $\lambda \in \R \backslash \Q$, then $Lf = 0$, by means of Fourier series, implies that $f$ is a constant, this means that $H^{0,0}_{}(\mathbb{T}; \mathfrak{v}) \cong \C$ and if $\lambda$ satisfies a Liouville condition, then $H^{0,1}_{}(\mathbb{T}; \mathfrak{v}) \cong \C$. 
    \end{exa}

    
        For the next example, we make use of the following basis for $\su(4)$.
        
        $$ \begin{matrix}
            T_1 = \begin{pmatrix}
                        i & 0 & 0 & 0 \\
                        0 & -i & 0 & 0 \\
                        0 & 0 & 0 & 0 \\
                        0 & 0 & 0 & 0 \\ 
                    \end{pmatrix}, ~ &
            T_2 = \begin{pmatrix}
                        i & 0 &   0 & 0 \\ 
                        0 & i &   0 & 0 \\ 
                        0 & 0 & -2i & 0 \\ 
                        0 & 0 &   0 & 0 \\ 
                    \end{pmatrix}, ~ &
            T_3 = \begin{pmatrix} 
                        i & 0 & 0 & 0 \\ 
                        0 & i & 0 & 0 \\ 
                        0 & 0 & i & 0 \\ 
                        0 & 0 & 0 & -3i \\ 
                        \end{pmatrix}, \\
            X_1 = \begin{pmatrix}
                        0 & i & 0 & 0 \\
                        i & 0 & 0 & 0 \\
                        0 & 0 & 0 & 0 \\
                        0 & 0 & 0 & 0 \\ 
                    \end{pmatrix}, ~ &
            X_2 = \begin{pmatrix}
                        0 & 0 & i & 0 \\
                        0 & 0 & 0 & 0 \\
                        i & 0 & 0 & 0 \\
                        0 & 0 & 0 & 0 \\
                    \end{pmatrix}, ~ &
            X_3 = \begin{pmatrix} 
                        0 & 0 & 0 & 0 \\
                        0 & 0 & i & 0 \\ 
                        0 & i & 0 & 0 \\ 
                        0 & 0 & 0 & 0 \\ 
                    \end{pmatrix}, \\
            Y_1 = \begin{pmatrix}
                        0 & -1 & 0 & 0 \\ 
                        1 &  0 & 0 & 0 \\ 
                        0 &  0 & 0 & 0 \\ 
                        0 &  0 & 0 & 0 \\ 
                    \end{pmatrix}, ~ &
            Y_2 = \begin{pmatrix} 
                        0 & 0 & -1 & 0 \\
                        0 & 0 &  0 & 0 \\ 
                        1 & 0 &  0 & 0 \\ 
                        0 & 0 &  0 & 0 \\ 
                    \end{pmatrix}, ~ &
            Y_3 = \begin{pmatrix}
                        0 & 0 &  0 & 0 \\ 
                        0 & 0 & -1 & 0 \\ 
                        0 & 1 &  0 & 0 \\ 
                        0 & 0 &  0 & 0 \\
                    \end{pmatrix}, \\
            X_4 = \begin{pmatrix}
                        0 & 0 & 0 & i \\
                        0 & 0 & 0 & 0 \\
                        0 & 0 & 0 & 0 \\
                        i & 0 & 0 & 0 \\ 
                    \end{pmatrix}, ~ &
            X_5 = \begin{pmatrix}
                        0 & 0 & 0 & 0 \\
                        0 & 0 & 0 & i \\
                        0 & 0 & 0 & 0 \\
                        0 & i & 0 & 0 \\
                    \end{pmatrix}, ~ &
            X_6 = \begin{pmatrix} 
                        0 & 0 & 0 & 0 \\
                        0 & 0 & 0 & 0 \\ 
                        0 & 0 & 0 & i \\ 
                        0 & 0 & i & 0 \\ 
                    \end{pmatrix}, \\
            Y_4 = \begin{pmatrix}
                        0 & 0 & 0 & -1 \\ 
                        0 & 0 & 0 &  0 \\ 
                        0 & 0 & 0 &  0 \\ 
                        1 & 0 & 0 &  0 \\ 
                    \end{pmatrix}, ~ &
            Y_5 = \begin{pmatrix} 
                        0 & 0 & 0 &  0 \\
                        0 & 0 & 0 & -1 \\ 
                        0 & 0 & 0 &  0 \\ 
                        0 & 1 & 0 &  0 \\ 
                    \end{pmatrix}, ~ &
            Y_6 = \begin{pmatrix}
                        0 & 0 & 0 &  0 \\ 
                        0 & 0 & 0 &  0 \\ 
                        0 & 0 & 0 &  -1 \\ 
                        0 & 0 & 1 &  0 \\
                    \end{pmatrix}. \\
        \end{matrix} $$

        Now we take $L_j = X_j - iY_j$ for $j = 1,2,3$ and, with an easy, but tedious, computation of commutators, we get that
        \begin{equation} \label{eq:eigenvectors}
            [T_j,L_k] = -i\lambda_{jk}L_k
        \end{equation}
        with $\lambda_{jk}$ real numbers and 
        
        \begin{tabular}{llllll}
            $\lambda_{1,1} = 2$, & $\lambda_{1,2} = 1$, & $\lambda_{1,3} = -1$, & $\lambda_{1,4} = 1$, & $\lambda_{1,5} = -1$, & $\lambda_{1,6} = 0$, \\
            $\lambda_{2,1} = 0$, & $\lambda_{2,2} = 3$, & $\lambda_{2,3} = 3$, & $\lambda_{2,4} = 1$, & $\lambda_{2,5} = 1$, & $\lambda_{2,6} = -2$, \\
            $\lambda_{3,1} = 0$, & $\lambda_{3,2} = 0$, & $\lambda_{3,3} = 0$, & $\lambda_{3,4} = 4$, & $\lambda_{3,5} = 4$, & $\lambda_{3,6} = 4$. \\
        \end{tabular}
        
        In order to compute the Levi form of subalgebras, we need to compute brackets between these $L_j$ and $\bar{L_k}$. This is done on Table \ref{tab:levi_form}.
        
        \begin{table}
            \caption{Each cell is the Lie bracket between the elements of the corresponding line and column.}
            \label{tab:levi_form}
            \begin{center} \begin{tabular}{c||c|c|c|c|c|c|}
                 & $\bar L_1$ & $\bar L_2$ & $\bar L_3$ & $\bar L_4$ & $\bar L_5$ & $\bar L_6$ \\
                 \hline
                 \hline
                 $ L_1 $ & $4iT_1$ & $-2iL_3$ & 0 & $2iL_5$ & 0 & 0 \\
                 $ L_2 $ & - & $2i(T_1 + T_2)$ & $2i \bar L_1$ & 0 & 0 & 0 \\
                 $ L_3 $ & - & - & $2i(T_2 - T_1$) & 0 & $-2i \bar L_6$ & 0 \\
                 $ L_4 $ & - & - & - & $\frac{2i}{3}(3T_1 + T_2 + 2 T_3)$ & $-2i\bar L_1$ & $-2i \bar L_2$ \\
                 $ L_5 $ & - & - & - & - & $\frac{2i}{3}(T_2 - 3T_1 + 2 T_3)$ & $-2i \bar L_3$ \\
                 $ L_6 $ & - & - & - & - & - & $\frac{4i}{3}(T_2 - T_3)$ \\
            \end{tabular} \end{center}
        \end{table}
        
        
        To simplify the computation, we denote by $\{\tau_1,\tau_2,\tau_3, \xi_1, \eta_1, \ldots, \xi_6, \eta_6 \}$ the dual basis of $\{T_1, T_2, T_3, X_1, Y_1 \ldots, X_6, Y_6\}$. We denote by $\mathbb T^d$ the $d$-dimensional torus. We write coordinates $(t_1, \ldots, t_d)$, the vector fields by $\del / \del t_1, \ldots, \del / \del t_d$ with dual basis denoted by $\theta_1. \ldots, \theta_d$.
        
        \begin{exa} \label{exe:cr_flat_not_a_product}
            Let  $G = \mathbb T^2 \times SU(4)$ and $Z_1 = T_1 - iT_2$ and $Z_2 = \del / \del t_1 + iT_3$. We define $\mathfrak h = \operatorname{span}_\C\{Z_1 Z_2,L_1, \ldots, L_3\}$. This structure is of type CR0.
        
            Now we will compute its Levi-form. If $\xi = \in \mathfrak h^\perp \cap \su(4)^*$, then $\xi = t(\theta_2)$ with $t \in \R \backslash \{0\}$, then $$ \mathcal L_\xi (X,Y) = \frac{1}{2i}\theta_2([Z,\bar Y])$$ and we easily see that it is flat. In fact, we just need to show that it vanishes on the basis. Clearly we have that $\theta_2([Z_1,\overline Z_2]) = \theta_2([Z_2,\overline Z_1]) = 0$ because $[Z_1, \overline Z_2] = [Z_2, \overline Z_1] = 0$. From Table \ref{tab:levi_form} we get that $\theta_2 (L_j, \overline L_k) = 0$ and finally from \eqref{eq:eigenvectors} we get that $\theta_2 (Z_j,\overline L_k) = 0$ because $[Z_j, \overline L_k] = \lambda \overline L_k$.
        \end{exa}

    \section{Invariant CR structures of maximal rank on the torus}
    \label{sec:cr_struct_torus}
    
    Let $\mathbb{T}$ be a torus of dimension $N = 2n + 1$. Assume it is endowed with a left-invariant CR structure $\mathfrak{m} \subset \C \mathfrak{t}$ of rank $n$. Since $\mathbb{T}$ is abelian, the structure $\mathfrak{m}$ is bi-invariant. We choose a basis $\{L_1, \ldots, L_n\}$ for $\mathfrak{m}$ which we complete to a basis for $\C \mathfrak{t}$ denoted by $\{L_1, \ldots, L_n, L_{n+1}, \ldots, L_{2n+1}\}$. Notice that we can take $L_{n+j} = \overline{L_j}$ for $j = 1, \ldots, n$ and take $L_{2n+1}$ to be real. We denote by $\{\tau_1. \ldots, \tau_{2n+1}\}$ the basis dual to $\{L_1, \ldots, L_n, L_{n+1}, \ldots, L_{2n+1}\}$. For a ordered multi-index $I$, we write $L_I = (L_{i_1}, \ldots, L_{i_q})$ and $\tau_I = \tau_{i_1} \wedge \ldots \wedge \tau_{i_q}$. 
    
    With the notation we just described, each $u \in C^\infty(\mathbb{T}, \Lambda^{0,q})$ can be written as $$u = \sum_{|I| = q} u_{I} \tau_I$$ and we have an expression for $\dbarb$ acting on $u$: $$ \dbarb u = \sum_{k = 1}^n \sum_{|I| = q} L_k u_{I} \tau_k \wedge \tau_I$$

We want to prove that under certain conditions on the CR structure $\mathfrak m$, for each cohomology class $[u] \in H^{0,q}_{\mathcal C}(\mathbb T; \mathfrak m)$ there exist a bi-invariant form $u_0$ such that $u - u_0 = \dbar_b v$ for some $v \in C^{\infty}(\mathbb T; \Lambda^{0,q-1})$, that is, we have that $[u] = [u_0]$ in $H^{0,q}_{}(\mathbb T; \mathfrak m)$.

Since we are assuming that each $L_j$ is bi-invariant, we can write it as $L_j = \sum_{j=1}^N a_{jk} \del / \del x_k$ with $a_{jk} \in \C$ and we denote write the symbol of $L_j$ as $$\widehat{L_j}(\xi) = i \sum_{j=1}^N a_{jk}\xi_k$$ for $\xi \in \Z^N$. Next, we state a condition that is sufficient for obtaining such bi-invariant representatives in every bi-degree.

\begin{dfn}
    \label{dfn:DC}
    We say that a CR subalgebra $\mathfrak{m} \subset \C \mathfrak{t}$ satisfies the {\bf(DC)} condition if there exist a basis $\{L_1, \ldots, L_n\}$ for $\mathfrak{m}$ and constants $C, \rho > 0$ such that
    $$ \max_j |\widehat{L_j}(\xi)| \geq C(1 + |\xi|)^{-\rho}, \quad \forall \xi \in \Z^N.$$
\end{dfn}

The definition of the {\bf(DC)}  condition does not depend on the choice of basis, in fact, let $\{L_1', \ldots, L_1'\}$ be any other basis for $\mathfrak{m}$, then there are constants $a_{jk} \in \C$ such that $L'_j = \sum_{k=1}^n a_{jk}L_k$ and we have
$$ |\widehat{L_j'}(\xi)| \leq \sum_{k=1}^n |a_{jk}|\widehat{L_k}(\xi)| \leq nA \max_j |\widehat{L_j}(\xi)|$$
with $A = \max_{j,k} |a_{jk}|$. That is, there are constants $c, C > 0$ such that
$$ c \max_j |\widehat{L_j'}(\xi)| \leq \max_j |\widehat{L_j}(\xi)| \leq C \max_j |\widehat{L_j'}(\xi)|$$

\subsection{Fourier series on forms}

Some of the computations in this section were inspired by \cite{dattori2016cohomology} and \cite{bergamasco1999global}.

For $u \in C^\infty(\mathbb{T}; \Lambda^{0,q})$, we write $u = \sum_{|I| = q} u_{I} \tau_I$ and since each $u_{I}$ is a smooth function on $\mathbb T$, by using Fourier transform, we can write $$ u(x) = \sum_{\xi \in \Z^N} e^{i \xi x} \sum_{|I| = q}  \widehat{u_{I}}(\xi) \tau_I = \sum_{\xi \in \Z^N} e^{i \xi x} \widehat{u}(\xi)$$ with $$\widehat{u}(\xi) = \sum_{|I| = q}  \widehat{u_{I}}(\xi) \tau_I.$$

We denote by $\|\cdot\|$ the value $$ \|\widehat{u}(\xi)\| = \max_{|J| = q} \{|\widehat{u_J}(\xi)|\}.$$ Since $u \in C^\infty(\mathbb{T}; \Lambda^{0,q})$, each $u_{I}$ is smooth, so for each $\nu \in \Z_+$ there is $C_\nu > 0$ such that $$ \|\widehat{u}(\xi)\| \leq C_\nu(1 + |\xi|)^{-\nu}$$ for all $\xi \in \Z^N.$
 
\begin{lem}
A $(0,q)$-form $u \in C^\infty(\mathbb{T}; \Lambda^{0,q)}$ is $\dbarb$-closed if, and only if,
    \begin{equation}
    \label{eq:fourier_closed}
        \left( \sum_{k=1}^n \widehat{L_k}(\xi) \tau_k\right) \wedge \widehat{u}(\xi) = 0,~\forall \xi \in \Z^N.
    \end{equation}
\end{lem}

\begin{proof}
    We have
    $$
    \begin{aligned}
        \dbarb u &= \sum_{|I| = q} \sum_{k = 1}^n \sum_{\xi \in \Z^N} \widehat{L_k u_{I}}(\xi)e^{i\xi x} \tau_k \wedge \tau_I \\
    & = \sum_{\xi \in \Z^N} \left( \sum_{k = 1}^n \sum_{|I| = q} \widehat{L_k u_{I}}(\xi) \tau_k \wedge \tau_I \right) e^{i \xi x} \\
    & = \sum_{\xi \in \Z^N} \left[  \left( \sum_{k = 1}^n \widehat{L_k}(\xi) \tau_k \right)\wedge \widehat{u}(\xi) \right] e^{i \xi x}
    \end{aligned}$$
    
    From the last equality we easily see that $\dbarb u = 0$ if and only if $\left( \sum_{k = 1}^n \widehat{L_k}(\xi) \tau_k \right)\wedge \widehat{u}(\xi) = 0$ for all $\xi \in \Z^N$.
\end{proof}

For $J = (j_1, \ldots, j_q)$ and $\sigma = j_k$, we write $J\backslash\sigma$ to denote the multi-index $(j_1, \ldots, \hat{j_k}, \ldots, j_q)$ and let $\varepsilon_{\sigma,J}$ be the signature of the permutation $(\sigma, J\backslash\sigma)$. 

\begin{lem}
    \label{lem:solvability_in_each_frequency}
    Let $u \in C^\infty(\mathbb{T}; \Lambda^{0,q})$ be a $\dbarb$-closed and suppose that $\sum_{k=1}^n \widehat{L_k}(\xi) \tau_k \neq 0$ and let $\sigma \in \{1, \ldots, n\}$ be such that $|\widehat{L_\sigma}(\xi)| = \max \{|\widehat{L_k}(\xi)|: k = 1, \ldots, n\}$. Then the $(q-1)$-form $$ \widehat{v}(\xi) = \sum_{|J| = q;~ \sigma \in J} \varepsilon_{\sigma,J} \frac{1}{\widehat{L_\sigma}(\xi)} \widehat{u_J}(\xi)\tau_{J \backslash \sigma}$$
    satisfies
    \begin{equation}
    \label{lem:dbarb_formal_solution}
            \left(\sum_{k=1}^n \widehat{L_k}(\xi) \tau_k\right) \wedge \widehat{v}(\xi) = \widehat{u}(\xi)
    \end{equation} 
    and there exist $C > 0$ not depending on $\xi$ such that
    \begin{equation}
    \label{ineq:dbarb_est_coef}
        \|\widehat{v}(\xi)\| \leq \frac{C}{|\widehat{L_\sigma}(\xi)|}\|\widehat{u}(\xi)\|.
    \end{equation}
\end{lem}


\begin{proof}
    Notice that Equation \eqref{ineq:dbarb_est_coef} follows directly from the definition of $\hat v(\xi)$.
    To prove Equation \eqref{lem:dbarb_formal_solution}, we adapt the proof of Lemma 2.1 of \cite{bergamasco1999global}. By definition, we have
    \begin{equation}
    \begin{aligned}
        \left(\sum_{k=1}^n \widehat{L_k}(\xi) \tau_k\right) \wedge \widehat{v}(\xi) &= \left(\sum_{k=1}^n \widehat{L_k}(\xi) \tau_k\right) \wedge \sum_{|J| = q;~ \sigma \in J} \varepsilon_{\sigma,J} \frac{1}{\widehat{L_\sigma}(\xi)} \widehat{u_J}(\xi)\tau_{J \backslash \sigma}. 
    \end{aligned}
    \end{equation}
    
    On the other hand, we can write $$\tau_\sigma = \frac{1}{\widehat{L}_\sigma (\xi)}\sum_{k=1}^n \widehat{L_k}(\xi) \tau_k -  \frac{1}{\widehat{L}_\sigma (\xi)}\sum_{k=1, k \neq \sigma}^n \widehat{L_k}(\xi) \tau_k$$ and now
    \begin{equation}
    \begin{aligned}
        \hat u(\xi) &= \sum_{|J|=q; \sigma \in J} \widehat{u_J}(\xi)\varepsilon_{\sigma,J} \tau_\sigma \wedge \tau_{J \backslash \sigma } + \sum_{|J|=q; \sigma \notin J} \widehat{u_J}(\xi) \tau_{J} \\
            &= \sum_{|J|=q; \sigma \in J} \widehat{u_J}(\xi)\varepsilon_{\sigma,J} \left(\frac{1}{\widehat{L}_\sigma (\xi)}\sum_{k=1}^n \widehat{L_k}(\xi) \tau_k\right) \wedge \tau_{J \backslash \sigma } \\
            &+ \sum_{|J|=q; \sigma \notin J} \widehat{u_J}(\xi) \tau_{J} - \sum_{|J|=q; \sigma \in J} \widehat{u_J}(\xi)\varepsilon_{\sigma,J} \left(\frac{1}{\widehat{L}_\sigma (\xi)}\sum_{k=1, k \neq \sigma}^n \widehat{L_k}(\xi) \tau_k\right) \wedge \tau_{J \backslash \sigma }
    \end{aligned}
    \end{equation}
    
    We define $$\hat u_\sigma(\xi) = \sum_{|J|=q; \sigma \in J} \widehat{u_J}(\xi)\varepsilon_{\sigma,J} \left(\frac{1}{\widehat{L}_\sigma (\xi)}\sum_{k=1}^n \widehat{L_k}(\xi) \tau_k\right) \wedge \tau_{J \backslash \sigma }$$ and $$\widehat{u}_{-\sigma}(\xi) = \sum_{|J|=q; \sigma \notin J} \widehat{u_J}(\xi) \tau_{J} - \sum_{|J|=q; \sigma \in J} \widehat{u_J}(\xi)\varepsilon_{\sigma,J} \left(\frac{1}{\widehat{L}_\sigma (\xi)}\sum_{k=1, k \neq \sigma}^n \widehat{L_k}(\xi) \tau_k\right) \wedge \tau_{J \backslash \sigma }.$$ Now we have $\hat u(\xi) = \hat u_{\sigma}(\xi) + \hat u_{-\sigma}(\xi)$. Since $u$ is $\dbar_b$-closed, it holds that $\left(\sum_{k=1}^n \widehat{L_k}(\xi) \tau_k\right) \wedge \hat u(\xi) = 0$ and, by construction, we have that $$ \left(\sum_{k=1}^n \widehat{L_k}(\xi) \tau_k\right) \wedge \hat u_{\sigma}(\xi) = 0$$ and thus $$ \left(\sum_{k=1}^n \widehat{L_k}(\xi) \tau_k\right) \wedge \hat u_{-\sigma}(\xi) = 0.$$ By expanding this last equality, and using the fact that $\tau_k \wedge \tau_J$ form a basis for all $(q+1)$-forms, we obtain that $\hat u_{-\sigma}(\xi) = 0$ and thus $\hat u(\xi) = \hat u_\sigma(\xi) = \left(\sum_{k=1}^n \widehat{L_k}(\xi) \tau_k\right) \wedge \widehat v(\xi)$.
\end{proof}

\begin{lem}
    Let $u \in C^\infty(\mathbb{T}; \Lambda^{0,q})$ $\dbarb$-closed and suppose that $\mathfrak{m}$ satisfies condition {\bf (DC)}. Then the form $$u_* = \sum_{\xi \in \Z^N \backslash \{0\}} e^{i\xi x} \widehat{u}(\xi)$$ is $\dbarb$-exact. 
\end{lem}

\begin{proof}
    For each $\xi \in \Z^N$ with $\xi \neq 0$, by the condition {\bf (DC)} we have $\sum_{k=1}^n \widehat{L_k}(\xi) \tau_k \neq 0$, and by Lemma \ref{lem:solvability_in_each_frequency}, there exist $\widehat v(\xi)$ such that $\left(\sum_{k=1}^n \widehat{L_k}(\xi) \tau_k\right) \wedge \widehat{v}(\xi) = \widehat{u}(\xi)$ and so we can define $v = \sum_{\xi \in \Z^N \backslash \{0\}} e^{i x \xi}$. Notice that by $\eqref{ineq:dbarb_est_coef}$ the form $v$ is smooth and by construction we have $\dbarb v = u_*$.
\end{proof}

If $[u]$ is a cohomology class in $H^{0,q}(T; \mathfrak m)$ we have that, if we define $u_0 = \widehat{u}(0) = \sum_{|I| = q} \widehat{u_J}(0) \tau_J$, then $[u] = [u_0]$. In fact, $u_* = u - u_0$ is $\dbarb$-exact by the last lemma and it follows that each cohomology class has a left-invariant representative. In particular, the cohomology groups $H^{0,q}(T; \mathfrak m)$ are finite-dimensional.

\section{Open problems}
\label{sec:op}
    Here we list a few open problems that can give us some directions for future research.
   
    \begin{enumerate}
        \item The main tool we used in this project was the Leray-Hirsch theorem which can be understood as a restricted version of the the Leray spectral sequence. What more could be obtained if the spectral sequence was used in full generality? 

        \item We imposed a strong condition to guarantee the existence of a cohomology extension. Is it possible to find weaker conditions that still guarantee the existence of cohomology extensions?
        
        \item The classification theorem for CR algebras of maximal rank was very useful in the study of the related cohomology spaces. Therefore, it is only natural to ask: is it possible to find a classification theorem for left-invariant CR structures in general, without assuming that they are of maximal rank?
    
        \item Let $G$ be a compact Lie group endowed with an elliptic Lie algebra $\mathfrak h$. We assume that there exist a maximal torus $T \subset G$ such that $\mathfrak h$ can be decomposed as \begin{equation} \label{eq:elliptic_toric} \mathfrak h = \mathfrak e \oplus \bigoplus_{\alpha \in \Delta_+} \mathfrak g_\alpha\end{equation} with $\mathfrak e$ a bi-invariant elliptic structure over $T$. Since every elliptic structure is Levi-flat, the smooth bundle $T \to G \to G / T$ admits local trivializations that are compatible with the elliptic structures. We also have that every cohomology class on the torus $T$ has a bi-invariant representative. Therefore, the proof of Theorem \ref{thm:main} can easily be adapted to this context. This raises the following question: what are the necessary and sufficient conditions so that left-invariant elliptic structures $\mathfrak h$ admit a decomposition as in \eqref{eq:elliptic_toric}? 
    \end{enumerate}

\bibliographystyle{alpha}
\bibliography{bibliography.bib}

\begin{thebibliography}{GMnMM13}

\bibitem[AI01]{alexandrov2001vanishing}
B.~Alexandrov and S.~Ivanov.
\newblock Vanishing theorems on {H}ermitian manifolds.
\newblock {\em Differential Geometry and its Applications}, 14(3):251--265,
  2001.

\bibitem[BCH08]{berhanu2008introduction}
S.~Berhanu, P.~D. Cordaro, and J.~Hounie.
\newblock {\em An introduction to involutive structures}, volume~6 of {\em New
  Mathematical Monographs}.
\newblock Cambridge University Press, Cambridge, 2008.

\bibitem[Bes08]{besse2008einstein}
A.~L. Besse.
\newblock {\em Einstein manifolds}.
\newblock Classics in Mathematics. Springer-Verlag, Berlin, 2008.
\newblock Reprint of the 1987 edition.

\bibitem[BJ19]{bor2019left-invariant}
G.~Bor and H.~Jacobowitz.
\newblock Left-invariant {CR} structures on 3-dimensional {L}ie groups.
\newblock {\em arXiv preprint arXiv:1909.08160}, 2019.

\bibitem[BP99]{bergamasco1999global}
A.~P. Bergamasco and G.~Petronilho.
\newblock Global solvability of a class of involutive systems.
\newblock {\em J. Math. Anal. Appl.}, 233(1):314--327, 1999.

\bibitem[CK04]{charbonnel2004classification}
J-Y. Charbonnel and H.~O. Khalgui.
\newblock Classification des structures {CR} invariantes pour les groupes de
  {L}ie compacts.
\newblock {\em J. Lie Theory}, 14(1):165--198, 2004.

\bibitem[DdSM16]{dattori2016cohomology}
P.~L. Dattori~da Silva and A.~Meziani.
\newblock Cohomology relative to a system of closed forms on the torus.
\newblock {\em Math. Nachr.}, 289(17-18):2147--2158, 2016.

\bibitem[GMnMM13]{gadea2013analysis}
P.~M. Gadea, J.~Mu\~{n}oz Masqu\'{e}, and I.~V. Mykytyuk.
\newblock {\em Analysis and algebra on differentiable manifolds}.
\newblock Problem Books in Mathematics. Springer, London, second edition, 2013.
\newblock A workbook for students and teachers, With a foreword by Andrew
  Swann.

\bibitem[Hal15]{hall2015lie}
B.~Hall.
\newblock {\em Lie groups, {L}ie algebras, and representations}, volume 222 of
  {\em Graduate Texts in Mathematics}.
\newblock Springer, Cham, second edition, 2015.
\newblock An elementary introduction.

\bibitem[HN12]{hilgert2011structure}
J.~Hilgert and K-H. Neeb.
\newblock {\em Structure and geometry of {L}ie groups}.
\newblock Springer Monographs in Mathematics. Springer, New York, 2012.

\bibitem[Pit88]{pittie1988dolbeault}
H.~V. Pittie.
\newblock The {D}olbeault-cohomology ring of a compact, even-dimensional {L}ie
  group.
\newblock {\em Proc. Indian Acad. Sci. Math. Sci.}, 98(2-3):117--152, 1988.

\bibitem[Spa95]{spanier1994algebraic}
E.~H. Spanier.
\newblock {\em Algebraic topology}.
\newblock Springer-Verlag, New York, [1995].
\newblock Corrected reprint of the 1966 original.

\bibitem[Tre92]{treves1992hypo}
F.~Treves.
\newblock {\em Hypo-analytic structures}, volume~40 of {\em Princeton
  Mathematical Series}.
\newblock Princeton University Press, Princeton, NJ, 1992.
\newblock Local theory.

\end{thebibliography}

\end{document}